\documentclass[smallextended]{svjour3}
\usepackage{amsfonts}
\usepackage{latexsym}
\usepackage{amsmath}
\usepackage{amssymb}
\usepackage{amscd}
\usepackage{epsfig}
\usepackage{multirow}
\usepackage{graphicx}
\addtocounter{MaxMatrixCols}{4}

\newcommand{\undertilde}[1]{\ensuremath{\mathord{\vtop{\ialign{##\crcr
   $\hfil\displaystyle{#1}\hfil$\crcr\noalign{\kern1.5pt\nointerlineskip}
   $\hfil\tilde{}\hfil$\crcr\noalign{\kern1.5pt}}}}}}
\begin{document}

\title{Deterministic ``Snakes and Ladders" Heuristic for the Hamiltonian Cycle Problem}
\author{P.~Baniasadi, V.~Ejov, J.A.~Filar, M.~Haythorpe, S.~Rossomakhine}

\institute{P. Baniasadi
\at Flinders University\\
\email{pouya.baniasadi@flinders.edu.au} \and V. Ejov
\at Flinders University\\
\email{vladimir.ejov@flinders.edu.au} \and J.A. Filar
\at Flinders University\\
\email{jerzy.filar@flinders.edu.au} \and M. Haythorpe (Corresponding author)
\at Flinders University, Ph: +61 8 8201 2834, Fax: +61 8 8201 2904\\
\email{michael.haythorpe@flinders.edu.au} \and  S. Rossomakhine
\at Flinders University\\
\email{serguei.rossomakhine@flinders.edu.au}} \maketitle
{\abstract
We present a polynomial complexity, deterministic, heuristic for solving the Hamiltonian Cycle Problem (HCP) in an undirected graph of order $n$. Although finding a Hamiltonian cycle is not theoretically guaranteed, we have observed that the heuristic is successful even in cases where such cycles are extremely rare, and it also performs very well on all HCP instances of large graphs listed on the TSPLIB web page.
The heuristic owes its name to a visualisation of its iterations. All vertices of the graph are placed on a given circle in some order. The graph's edges are classified as either snakes or ladders, with snakes forming arcs of the circle and ladders forming its chords.  The heuristic strives to  place exactly $n$ snakes on the circle, thereby forming a Hamiltonian cycle.
The Snakes and Ladders Heuristic (SLH) uses transformations inspired by $k-$opt algorithms such as the, now classical, Lin--Kernighan heuristic to reorder the vertices on
the circle in order to transform some ladders into snakes and {\em vice versa}. The use of a suitable stopping criterion ensures the heuristic terminates in
polynomial time if no improvement is made in $n^3$ major iterations.
\keywords{Hamiltonian cycle problem, Lin-Kernighan, heuristic, polynomial time, connected graph}}

\section{Introduction}\label{sec-Introduction}

The essence of the Hamiltonian Cycle Problem (HCP, for short) is
contained in the following, deceptively simple statement: {\em Given a graph $G$, find a Hamiltonian cycle (HC) (a simple cycle that contains
all vertices of the graph) or prove that such a cycle does not exist.} A graph is said to be {\em Hamiltonian}\footnote{The name stems from Sir William Hamilton's
investigations of such cycles on the dodecahedron graph around
1856, however Leonhard Euler studied the famous ``knight's tour" on a
chessboard as early as 1759.} if it contains at least one Hamiltonian cycle, and non-Hamiltonian otherwise.

Henceforth, a  {\em graph\/} of order $n$ will mean an $n$-vertex graph which is both simple (without self loops or multiple edges) and undirected (every edge admits two-way traffic.)

The HCP is known to be NP-complete and has become a challenge that
attracts mathematical minds both in its own right and because of
its close relationship to the famous {\em Travelling Salesman
Problem\/} (TSP). An efficient solution of the latter would have
an enormous impact in operations research, optimisation and
computer science.  However, the TSP is merely the problem of identifying a Hamiltonian cycle (also called a ``tour") of minimal length, where each edge has an associated length (or weight), and the length of the cycle is the sum of the lengths of edges comprising that cycle.    The latter constitutes a simple, linear, objective function and it can be argued that much of the difficulty of the TSP is embedded in the HCP, namely, in finding an optimal tour in the space of Hamiltonian cycles.

There is now an extensive body of literature devoted to TSP and its many variants. Besides theoretical developments there are also many exact algortihms, and heuristics. The reader is referred to the comprehensive books by Lawler {\em et al} \cite{Lawl85Traveling} and Gutin and Punnen \cite{Gutin}. For many researchers, solution of HCP becomes a simple corollary of the TSP although, in principle, the former need not require the minimisation of any objective function.

Some of the most successful heuristics for solving TSP are local search methods that exploit the so-called ``$k-$opt" transformations that facilitate movement from one tour to a shorter tour via an exchange of exactly $k$ edges. The notion of a ``$2-$opt" transformation is widely attributed to Flood \cite{Flood}.   Subsequently, Lin \cite{Lin} and Lin and Kernighan (L--K) \cite{LK} developed powerful heuristics that exploited ``$k-$opt" transformations more extensively. Indeed,  L--K is still embedded in some of the best modern heuristics for the TSP, notably Helsgaun's LKH \cite{Helsgaun}.  For a comprehensive discussion of modern developments we refer the reader to Applegate {\em et al} \cite{ConcordeBook}.

The Snakes and Ladders Heuristic (SLH) for the Hamiltonian Cycle Problem presented here is a polynomial-complexity algorithm inspired by, but distinctly different from, the $k-$opt heuristics.  The name Snakes and Ladders Heuristic comes from our visual representation of iterations of the algorithm described in detail in the next section.

In order to attempt to find a Hamiltonian cycle of a given graph $G$, we place its vertices in some order on a circle. Any edges between adjacent vertices on the circle are rendered as arcs of the circle (which we call {\it snakes}), while all other edges are rendered as chords of the circle (which we call {\it ladders}). If two adjacent vertices on the circle have no snake between them, we say there is a {\em gap} between them. SLH attempts to place all edges of a Hamiltonian cycle on the circle by a number of transformations that are isomorphisms of the underlying graph. While some, but not all, of these isomorphisms coincide with $k-$opt transformations, the main difference compared to $k-$opt search methods is that our approach does not require an improvement of a TSP-type objective function; rather, SLH seeks changes in the arrangement of vertices of the graph on the circle with the goal of facilitating the eventual {\em closure} of gaps.  We might say that in SLH, $k-$opt transformations have been generalised to $k-$change transformations in an appropriate way.

Whereas a $k-$opt method attempts to ``improve" after each iteration by  decreasing the number of gaps until a Hamiltonian cycle is obtained; SLH  uses a more general notion of ``improvement" by trying to achieve a suitable balance between  increasing the number of Hamiltonian cycle edges on the circle and decreasing the number of gaps.  As a consequence, SLH differs from standard $k-$opt heuristics in three important ways.

\begin{enumerate}
\item[(1)] SLH performs a sequence of compositions of two simple generator operations: {\it gamma} ($\gamma$) and {\it kappa} ($\varkappa$). The sequence so constructed often results in transformations that - under $k-$opt heuristics - would not be allowed, or would be difficult to identify. We later prove that all Hamiltonian cycles in any given graph of size $n$ are reachable through the use of these two operations, with the number of operations required bounded above by a linear function of $n$.
\item[(2)] The $k-$opt heuristics update the tour only when an \lq\lq improvement" is found, while SLH allows {\em floating} and {\em opening} transformations that may result in either no improvement or a \lq\lq sacrifice", respectively.
\item[(3)] The $k-$opt heuristics rely on randomisation techniques to obtain a Hamiltonian cycle while SLH does not take advantage of these techniques, and is designed to run on any initial input arrangement in a deterministic fashion.
\end{enumerate}

As indicated in item 3, we implement SLH as a {\em deterministic heuristic}. That is, for a graph and a given starting orientation, the heuristic will produce the same output every time it is run. A stopping condition is chosen to ensure that SLH will terminate in polynomial time, either by identifying a Hamiltonian cycle, or failing to improve after $n^3$ iterations. SLH has been implemented in C++. Although SLH is not guaranteed to find a Hamiltonian cycle in a Hamiltonian graph, preliminary experiments on many graphs (not exceeding 5000 vertices) have succeeded in all cases. That is, a Hamiltonian cycle has been found in all graphs that were Hamiltonian, while termination declaring the graph to be ``likely non-Hamiltonian" was reached in all instances of graphs known, a priori, to not possess any Hamiltonian cycles.

This paper is organised as follows: Section \ref{sec-algorithm} introduces the basic idea behind our approach, the transformations of SLH and our algorithm implementation. Section \ref{sec-motivation} gives a plausible explanation of why SLH is effective in finding Hamiltonian cycles. Section \ref{sec-performance} reports on some of the experiments performed with the algorithm as well as a comparison with well known TSP solvers. Finally, our conclusions of this paper are presented in Section \ref{sec-conclusion}, including a link to our website, http://fhcp.edu.au/slhweb/ \cite{SLH}, where readers are invited to test SLH on either built-in, or user-supplied, problems.

\section{Description of the algorithm} \label{sec-algorithm}

We start by introducing the terminology that we use in this paper. Then we will describe the transformations that are used in our implementation of SLH. Finally, we discuss our implementation of the SLH algorithm.

\subsection{Basic idea}

In our approach, we place the vertices of a simple undirected graph $G$ on a circle in some order. Then, all edges of $G$ between adjacent vertices on the circle are represented as arcs on the circle, which we call \emph{snakes}, while the other edges are represented as chords of the circle, which we call \emph{ladders}.

The arrangements of vertices $\{1,\ldots,n \}$ on the circle form natural equivalence classes.
Namely, two arrangements are said to be equivalent if either one can be transformed to the other via a rotation
(clockwise or anti-clockwise), a reversal of the ordering, or a composition of both. This implies two arrangements are equivalent if and only if all vertices have the same neighbours in both arrangements. For example, the arrangement (1,3,4,6,2,5) is equivalent to $(2,5,1,3,4,6)$ and $(4,3,1,5,2,6)$ but not $(3,6,4,2,1,5)$.

It is clear that any member of the equivalence class, for a fixed graph $G$, contains the same snakes and ladders as any other member of this class. We use the term {\em ordering}, or {\em circle ordering} to denote such an equivalence class, and give it the symbol $\mathcal{C}_G$. However, since the algorithm and all definitions in this paper are given for a fixed graph, we henceforth drop the subscript $G$. For a given ordering $\mathcal{C}$, if two vertices $x$ and $y$ that are adjacent on the circle are not connected by a snake we say that $(x,y)$ is a \emph{gap} on $\mathcal{C}$. When no confusion is possible, we will also use the term ordering to denote a particular member of the equivalence class, and use special notation to describe the ordering. For any adjacent pair of vertices, there are three distinct possibilities; there may be a snake between them, a gap between them, or there may be either a gap or a snake. These situations will be denoted by $a\smallfrown b$, $a~|~b$, and $a, b$ respectively. This notation will allow us to define transformations that require the presence of particular snakes or gaps in certain parts of the ordering, but remain defined regardless of what is present elsewhere in the ordering. One example of this notation is $(a\smallfrown b,\ldots,c~|~d,\ldots,e,a)$, where $a$ is viewed as the initial vertex, and $e$ as the final vertex. In this example, $b$ immediately follows $a$, then some number (possibly zero) additional vertices follow $b$ before $c$ is reached. Likewise, $d$ immediately follows $c$, and then some number (possibly zero) additional vertices follow before reaching $e$, and then $a$ follows as the circle is closed. When denoting an ordering, the initial vertex is repeated at the end. Equivalent notation without
the initial vertex repeated refers to a {\em segment} of the ordering. For the sake of clarity, we can choose to view a given ordering as an ordered set of segments. For example, if we have an ordering $(a,\ldots,b,c,\ldots,d,a)$, we can define segments $A = (a,\ldots,b)$ and $B = (c,\ldots,d)$ and rewrite the ordering as $(A,B)$. We additionally define segment $A^R$ to be the reverse of $A$, so in the previous example $A^R = (b,\ldots,a)$ and ordering $(A^R,B) = (b,\ldots,a,c,\ldots,d,b)$.


We note that a Hamiltonian cycle corresponds to an ordering containing exactly $n$ snakes, or equivalently, an ordering with $0$ gaps.
Hence, the ultimate goal for SLH is to close all the gaps on $\mathcal{C}$ through a series of transformations. These transformations
can be described as compositions of two {\em generator transformations}, specifically the following two isomorphisms:
\begin{enumerate}
\item $\gamma$ isomorphism. This transformation, denoted by $\gamma(y,x,a)$, maps an ordering $(x,\ldots,b,a,\ldots,y,x)$
to ordering  $(b,\ldots,x,a,\ldots,y,b)$. That is, if we define segments $A = (x,\ldots,b)$ and $B = (a,\ldots,y)$, then $\gamma(y,x,a)$ maps ordering $(A,B)$ to ordering $(A^R,B)$. The definition implies that
$\gamma(y,x,a)$ is only defined for orderings in which $x$ and $y$ are adjacent on the ordering. Vertex $b$ is implicitly defined as being the vertex adjacent to $a$ on the segment $(y,x,\ldots,a)$.

\begin{figure}[h!]\begin{center}\includegraphics[scale=0.265]{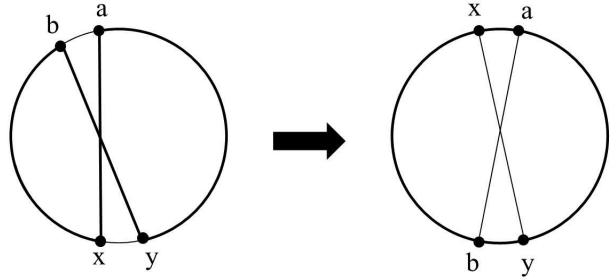}\\
\caption{Generator transformation $\gamma(y,x,a).$}\label{fig: op16}\end{center}\end{figure}

\item $\varkappa$ isomorphism. This transformation, denoted by $\varkappa(x,a,c,d)$, maps an ordering
\linebreak$(x,\ldots,e,c,\ldots,a,b,\ldots,f,d,\ldots,y,x)$ to an ordering $(e,\ldots,x,a,\ldots,c,d,\ldots,y,f,\ldots,b,e)$. That is, if we define segments $A = (x,\ldots,e)$, $B = (c\ldots,a)$, $C = (b,\ldots,f)$ and $D = (d,\ldots,y)$, then $\varkappa(x,a,c,d)$ maps ordering $(A,B,C,D)$ to ordering $(A^R,B^R,D,C^R)$. The definition implies that
$\varkappa(x,a,c,d)$ is only defined for orderings in which the segment $(x,\ldots,c,\ldots,a)$ does not contain $d$. This definition allows for
$e = x$, that is, segment $A$ contains a single vertex; for $y = d$, that is, segment $D$ contains a single vertex; and for segment $C$ to be empty (i.e. so vertex $a$ directly precedes vertex $d$ on the ordering), or for $b = f$, that is, segment $C$ contains a single vertex.

\begin{figure}[h!]\begin{center}\includegraphics[scale=0.20]{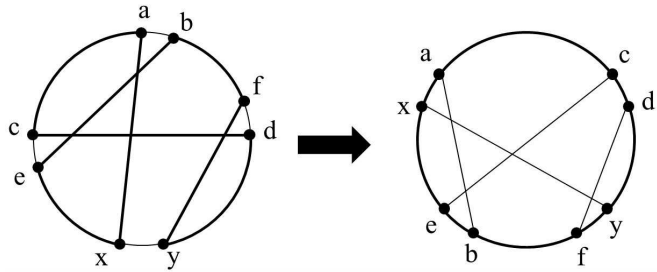}\\
\caption{Generator transformation $\varkappa(x,a,c,d)$.}\label{fig: op01}\end{center}\end{figure}
\end{enumerate}

Note that, although the figures for the generator transformations include ladders and snakes, these need not be present for the isomorphisms to be performed.

Later in Section \ref{sec-motivation}, we will show that for any initial ordering $\mathcal{C}$
and any Hamiltonian ordering $\mathcal{C}^H$, that is, an ordering corresponding to a Hamiltonian cycle $H$, there exists a transformation $T$
from $\mathcal{C}$ to $\mathcal{C}^H$ that is a composition of a number of $\gamma$ and $\varkappa$ transformations in some order. Namely,
\begin{equation*}\label{comp_gk}
T=\prod_{i=1}^\ell\gamma^{\varepsilon_i} \circ \varkappa^{\delta_i},
\end{equation*}
where $\varepsilon_i$ and $\delta_i$ can take values $0$ or $1$, and $\ell$ bounds the number of such transformations.
In fact, our algorithm attempts to iteratively build such a transformation $T$ to
some unknown Hamiltonian ordering as a composition of SLH transforms that are
themselves specially designed combinations of $\gamma$ and $\varkappa$ transformations. To ensure that
our algorithm is polynomially bounded we choose $\ell \leq n^4$.

Note that $T$ transforms the initial ordering to an ordering with zero gaps.
Therefore, SLH transformations that reduce the number
of gaps are viewed as desirable. Such transformations are widely used in TSP algorithms and are called
$k-$opt exchanges. Adapting $k-$opt transformations for HCP, they improve
the tour by exchanging $k$ ladders with gaps and snakes ($k$ in total) in such a way that
the new tour contains fewer gaps. Even though SLH incorporates $k-$opt
transformations, it also incorporates transformations that may preserve the weight of the tour (i.e. keep
the same number of gaps) or even increase the weight of the tour (i.e. increase
the number of gaps). We call these $k-$flo transformations, to emphasise their indeterminant (floating)
nature. Unlike $k-$opt transformations, $k-$flo transformations are defined for situations where it is not known,
in advance, whether certain gaps or ladders are present. Specifically, a $k-$flo transformation is one in which as many as $k$ ladders
may be turned into snakes, if they exist. Then, if it so happens that some of these ladders and/or gaps exist,
the transformation might produce a new ordering with a reduced number of gaps. In such a case, the $k-$flo transformation specialises to an
$m-$opt transformation for some $m \leq k$.

\subsection{Operations}

In our implementation, we found it useful to consider special combinations of $\gamma$ and $\varkappa$ that we found to be the most effective among the many alternatives we tried. These can be partitioned into three types; \emph{closing}, \emph{floating} and \emph{opening}. In each case, certain ladders, snakes and gaps are required to either be present or absent for the combination of generator transformations to produce the transformation.

{\bf SLH--closing transformations:} These correspond to particular 2--opt and 3--opt transformations:
\begin{enumerate}
\item SLH closing $2$--opt type 1 transformation is $\gamma(y,x,a)$. It maps ordering
$(x,\ldots,b~|~a,\ldots,y~|~x)$, where $(x,a)$ is a ladder, to an ordering $(b,\ldots,x\smallfrown a,\ldots,y~|~b)$.

\begin{figure}[h!]\begin{center}\includegraphics[scale=0.28]{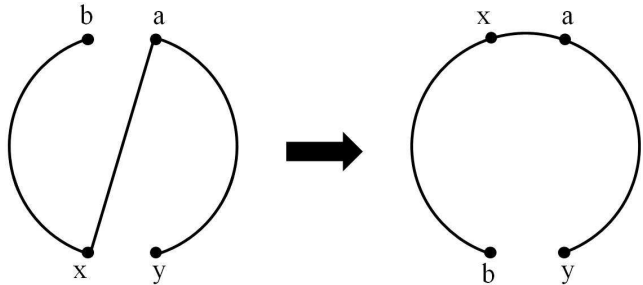}\\
\caption{SLH closing $2$--opt type 1 transformation, $\gamma(y,x,a)$.}\label{fig: op01}\end{center}\end{figure}

\item SLH closing $2$--opt type 2 transformation is $\gamma(y,x,a)$ that maps an ordering $(x,\ldots,b,a,\ldots,y~|~x)$,
where $(x,a)$ and $(y,b)$ are ladders, to an ordering $(b,\ldots,x\smallfrown a,\ldots,y\smallfrown b)$.

\begin{figure}[h!]\begin{center}\includegraphics[scale=0.21]{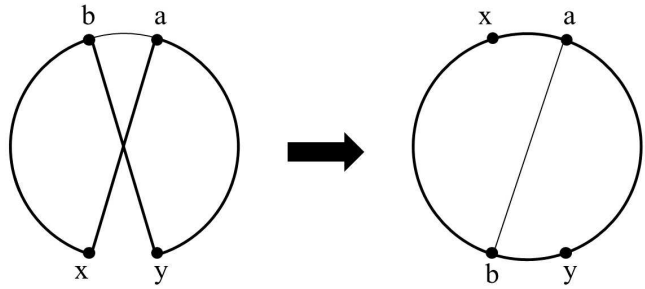}\\
\caption{SLH closing $2$--opt type 2 transformation, $\gamma(y,x,a)$.}\label{fig: op01}\end{center}\end{figure}

\item SLH closing $3$--opt transformation is $\gamma(c,y,b) \circ \gamma(y,x,a)$. It maps ordering \linebreak$(x,\ldots,c,a,\ldots,b,d,\ldots,y~|~x)$, where $(x,a)$, $(c,d)$ and $(b,y)$ are ladders, to an ordering \linebreak$(d,\ldots,y\smallfrown b,\ldots,a\smallfrown x,\ldots,c\smallfrown d)$.

\begin{figure}[h!]\begin{center}\includegraphics[scale=0.28]{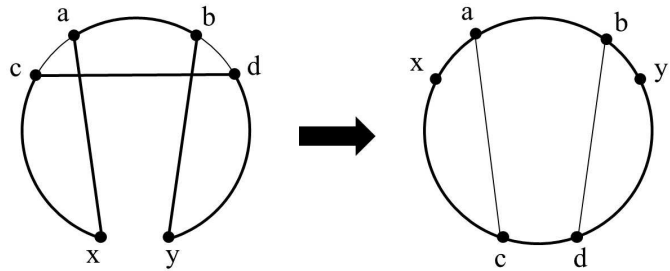}\\
\caption{SLH closing 3-opt transformation, $\gamma(c,y,b) \circ \gamma(y,x,a)$.}\label{fig: op02}\end{center}\end{figure}
\end{enumerate}

{\bf SLH--floating transformations:} These are special ladder--snake interchange transformations that generate new orderings,
where the number of gaps is either unchanged, or possibly reduced (if certain ladders, gaps, or both are present in the first ordering).
Ladders which (if they exist) will be turned into snakes after the transformation are represented by dashed lines in the first ordering,
with the corresponding snake also represented by a dashed line in the second ordering. Of course, if such a ladder does not exist
then the corresponding dashed line in the second ordering represents a gap. In particular:
\begin{enumerate}

\item SLH floating $2$--flo transformation is $\gamma(y,x,a)$ that maps an ordering $(x,\ldots,b,a,\ldots,y~|~x)$,
where $(x,a)$ is a ladder, to an ordering $(b,\ldots,x\smallfrown a,\ldots,y,b)$. If $(b,y)$ is also a ladder, this transformation
closes at least one gap.

\begin{figure}[h!]\begin{center}\includegraphics[scale=0.28]{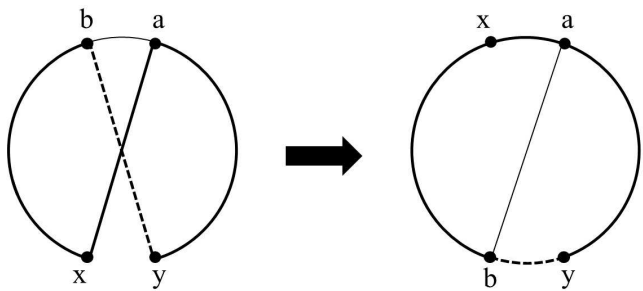}\\
\caption{SLH floating $2$--flo transformation, $\gamma(y,x,a)$.}\label{fig: op16}\end{center}\end{figure}

\item SLH floating $3$--flo transformation is $\gamma(c,y,b) \circ \gamma(y,x,a)$. It maps ordering \linebreak$(x,\ldots,c,a,\ldots,b,d,\ldots,y~|~x)$, where $(x,a)$ and at least
one of $(c,d)$ and $(b,y)$ are ladders, to an ordering $(d,\ldots,y,b,\ldots,a\smallfrown x,\ldots,c,d)$. If both $(c,d)$ and $(b,y)$ are ladders, this transformation closes at least one gap.

\begin{figure}[h!]\begin{center}\includegraphics[scale=0.21]{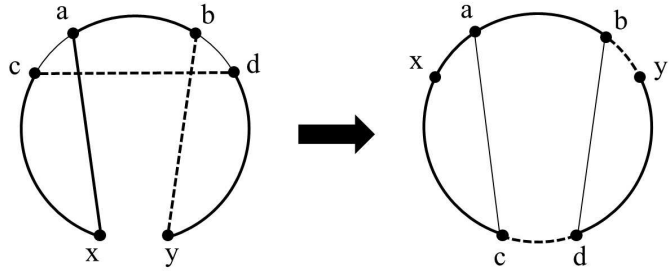}\\
\caption{SLH floating $3$--flo transformation, $\gamma(c,y,b) \circ \gamma(y,x,a)$.}\label{fig: op17}\end{center}\end{figure}

\item SLH floating $4$--flo type 1 transformation is $\gamma(d,b,y) \circ \gamma(b,e,f) \circ \varkappa(x,a,c,d)$.
It maps ordering  \linebreak$(x,\ldots,e,c,\ldots,a,b\ldots,d,f,\ldots,y~|~x)$, where $(x,a)$, $(b,y)$ and $(c,d)$ are ladders,
to an ordering \linebreak$(d\ldots,b,y,\ldots,f,e,\ldots,x\smallfrown a,\ldots,c,d)$. If $(e,f)$ is also a ladder, this transformation closes at least one gap.

\begin{figure}[h!]\begin{center}\includegraphics[scale=0.21]{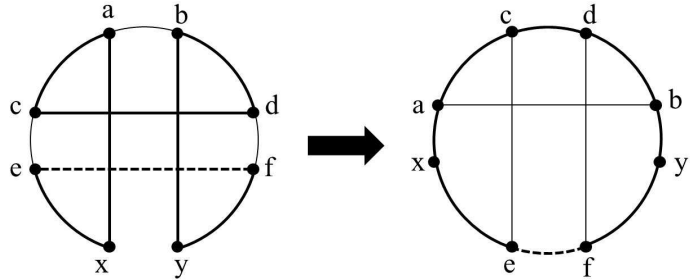}\\
\caption{SLH floating $4$--flo type 1 transformation,  $\gamma(d,b,y) \circ \gamma(b,e,f) \circ \varkappa(x,a,c,d)$.}\label{fig: op12}\end{center}\end{figure}

\item SLH floating $4$--flo type 2 transformation is $\varkappa(x,a,c,d)$ that maps ordering
\linebreak$(x,\ldots,e,c,\ldots,a,b\ldots,f,d,\ldots,y|x)$, where $(x,a)$ and $(c,d)$ are ladders, and at least
one of $(e,b)$ and $(f,y)$ are ladders, to an ordering $(f,\ldots,b,e,\ldots,x\smallfrown a\ldots,c,d,\ldots,y,f).$ If both $(e,b)$ and $(f,y)$ are ladders,
this transformation closes at least one gap.

\begin{figure}[h!]\begin{center}\includegraphics[scale=0.21]{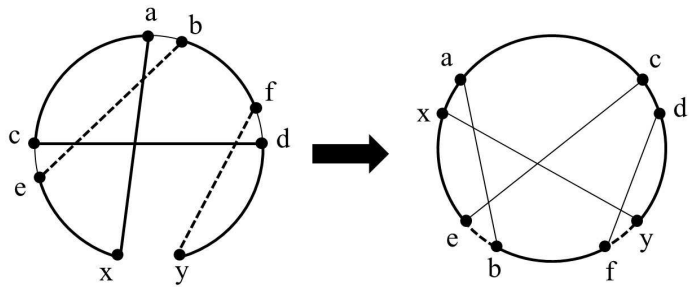}\\
\caption{SLH floating $4$--flo type 2 transformation, $\varkappa(x,a,c,d)$.}\label{fig: op15}\end{center}\end{figure}

\item  SLH floating $5$--flo transformation is $\varkappa(j,d,g,h) \circ \varkappa(x,a,c,d)$ that maps ordering \linebreak$(x,\ldots,e,c,\ldots,a,f,\ldots,g,b,\ldots,j,d,h\ldots,y~|~x)$,  where $(x,a)$, $(b,y)$, $(c,d)$ and $(f,e)$ are ladders, to an ordering
 $(b,\ldots,j,d,c,\ldots,a\smallfrown x,\ldots,e,f,\ldots,g,h,\ldots,y,b)$. If $(g,h)$ is also a ladder, this transformation closes at least one gap.

\begin{figure}[h!]\begin{center}\includegraphics[scale=0.21]{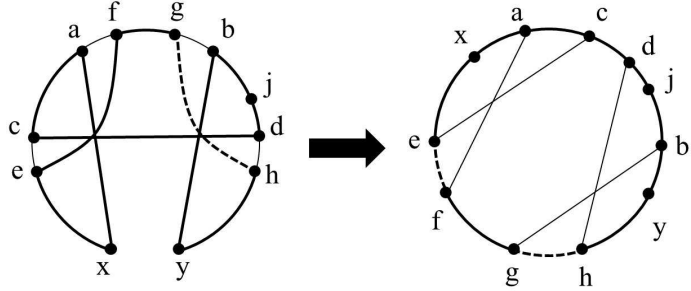}\\
\caption{SLH floating $5$--flo transformation, $\varkappa(j,d,g,h) \circ \varkappa(x,a,c,d)$.}\label{fig: op14}\end{center}\end{figure}
\end{enumerate}

{\bf SLH opening transformations:} Importantly, we introduce a single snake-ladder interchange transformation that generates a new ordering that generally contains one more gap.
\begin{enumerate}
\item SLH opening $4-$flo transformation is $\varkappa(x,a,c,d)$.
It maps ordering \linebreak$(x,\ldots,e,c,\ldots,a,b,\ldots,f,d,\ldots,y~|~x)$, where $(x,a)$ and $(c,d)$ are ladders, to an ordering \linebreak$(e,\ldots,x\smallfrown a,\ldots,c,d,\ldots,y,f,\ldots,b,e)$.
In general, this transformation increases the number of gaps by 1, though this is not always the case.

\begin{figure}[h!]\begin{center}\includegraphics[scale=0.21]{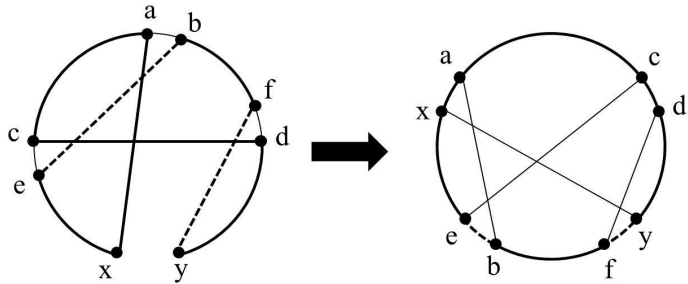}\\
\caption{SLH opening $4-$flo transformation $\varkappa(x,a,c,d)$.}\label{fig: op21}\end{center}\end{figure}
\end{enumerate}

Note that, in all cases, these transformations are designed to be performed for an ordering containing a gap $(x,y)$. To represent this, we say that a transformation is performed \lq\lq around gap $(x,y)$".

\subsection{Algorithm}
SLH works in four mains stages. In each stage it attempts to find a Hamiltonian cycle by a particular approach. If no Hamiltonian cycle is found within the (polynomial complexity) constraints of that stage, SLH moves on to the next stage in order to attempt a different approach. Finally, if stage 3 has failed to reduce the number of gaps after $n^3$ iterations, the graph is declared to be ``likely non-Hamiltonian". The following is a description of our algorithm implementation.

\vspace*{1cm}{\bf Stage 0.}
\begin{enumerate}
\vspace*{-0.4cm}\item[{\bf (0.1)}] Let the initial ordering be the original assignment of the graph.
\item[{\bf (0.2)}] Perform SLH-closing transformations to obtain new orderings. Continue until the number of gaps can not be reduced by any SLH-closing transformation. If the number of gaps is zero, stop: a Hamiltonian cycle has been found.
\item[{\bf (0.3)}] Create a gap list, initially empty, and an ordering list, initially containing only the most recently obtained ordering. These two lists will act as tabu lists. Go to stage 1.
\end{enumerate}

We note that we can substitute stage 0 with any effective polynomially bounded $k-$opt algorithm\footnote{An effective implementation of L--K, for example Helsgaun's LKH \cite{Helsgaun}, could be such an algorithm.}, disallowing the use of any randomisation techniques to retain the deterministic nature of SLH. The final ordering provided by such an algorithm will become the initial ordering for stage 1. Using operations from a sufficiently efficient $k-$opt algorithm instead of our closing operations may improve the performance of SLH. However, for this proof of concept, and to emphasise the effectiveness of SLH at stages 1, 2 and 3, we chose the simplest possible closing transformations.

{\bf Stage 1.}
\begin{enumerate}\vspace*{-0.4cm}\item[{\bf (1.1)}] Perform any SLH-floating transformation around a gap $g$ on the latest ordering. Save $g$ in the gap list and continue as follows:
\begin{enumerate}
\item[(a)] If the new ordering contains a gap that is not in the gap list, add the ordering to the ordering list and proceed to stage 1.2.
\item[(b)] If the new ordering contains only gaps that are in the gap list repeat stage 1.1 with a different SLH-floating transformation on the previous ordering. If no SLH-floating transformation produces an ordering which contains only gaps that are in the gap list already, then revert to an earlier ordering, and continue (perhaps needing to revert back to an even earlier ordering) until a new SLH-floating transformation can be performed that does produce an ordering containing a gap that is not in the gap list. Then, add the ordering to the ordering list and proceed to stage 1.2.
\item[(c)] If no SLH-floating transformation on any ordering in the ordering list produces an ordering with no gaps in the gap list, proceed to stage 2.
\end{enumerate}
\item[{\bf (1.2)}]If the number of gaps is less than in the previous iteration, then remove all orderings from the ordering list, and all gaps from the gap list. If a Hamiltonian cycle is found, stop.
\item[{\bf (1.3)}]Return to stage 1.1.
\end{enumerate}

Note that, according to the stopping rule,
if an SLH-floating transformation creates a gap that
has already been obtained by some previous transformation, the transformation will not be
performed and the ordering will not be listed.
Thus, the stopping rule limits the number of orderings in the ordering list at this stage
to $n^2.$ Once each ordering in the ordering list has been considered, we proceed to stage 2.

{\bf Stage 2.}
\begin{enumerate}
\vspace*{-0.4cm}\item[{\bf (2.1)}]Record the most recently obtained ordering $\mathcal{C}$ and the number of gaps $g(\mathcal{C})$.
\item[{\bf (2.2)}]Perform an SLH-opening transformation from a gap $g$ in $\mathcal{C}$ and save the ordering so obtained in the ordering list.
\item[{\bf (2.3)}]Repeat stage 1 until no gap can be closed. If the number of gaps is zero, stop: a Hamiltonian cycle has been found.
\item[{\bf (2.4)}]If the number of gaps in the most recent ordering is less than $g(\mathcal{C})$, empty the gap list, remove all orderings except the most recent ordering from the ordering list, and return to stage 1.
\item[{\bf (2.5)}]Recover ordering $\mathcal{C}$ and return to stage 2.2 but perform a different SLH-opening transformation. If all possible SLH-opening transformations from $g$ have been previously considered, go to stage 3.
\end{enumerate}

{\bf Stage 3.}
\begin{enumerate}
\vspace*{-0.4cm}\item[{\bf (3.1)}] Record the most recently obtained ordering $\mathcal{C}'$.
\item[{\bf (3.2)}] Perform an opening transformation from any gap $g$ in $\mathcal{C}'$ and add the resulting ordering to the ordering list. Then, attempt to perform $k-$opt transformations (floating transformations that reduce the number of gaps), without obtaining already listed orderings, saving each new ordering in the ordering list. Continue until no more gaps can be closed. If the number of gaps is zero, stop: a Hamiltonian cycle has been found. If at any stage
the number of orderings in the ordering list exceeds $n^3$, go to stage 3.5.
\item[{\bf (3.3)}] If the number of gaps in the most recently obtained ordering is less than $g(\mathcal{C})$ (recorded in stage 2.1), empty the gap list, remove all orderings except the most recent ordering from the ordering list, and return to stage 1. Otherwise, revert to the most recently obtained ordering, descending from $\mathcal{C}'$, in which $k-$opt transformations that have not yet been tried are possible. Continue attempting to perform $k-$opt transformations without obtaining already listed orderings, saving each new ordering in the ordering list.
\item[{\bf (3.4)}] If all possible $k-$opt transformations in each previous ordering descending from $\mathcal{C}'$ have already been tried, return to stage 3.1.
\item[{\bf (3.5)}] If the number of orderings in the ordering list is equal to $n^3$, and the number of gaps is greater than zero in all, declare the graph to be ``likely non-Hamiltonian" and stop.
\end{enumerate}

Note that in stages above where a transformation is selected, there may be multiple eligible transformations. The algorithm searches for any applicable transformations, and as soon as an eligible transformation is discovered, it is performed. If the transformations can be performed from any gap, the gaps are perused in order of the ordering. This ensures that the above process is deterministic, with the order of transformations determined entirely by the initial assignment of the graph.

\subsection{Worst-case algorithmic complexity}

In the algorithm implementation above, stage 3 has the largest search space and the stopping condition with the largest bound, so in the worst case dominates the execution time. In each iteration of stage 3 we allow as many opening transformations as are necessary to enable a floating transformation to close a gap. Then, there will be a sequence of $O(n)$ such floating transformations that close gaps. If $k$ is the maximum degree of the graph, then there are $O(nk^4)$ potential floating transformations at each step of the sequence. This number of transformations arises because the SLH floating 5-flo transformation requires four edges to be chosen, the first emanating from any vertex next to a gap, and each subsequent edge emanating from a vertex determined by the previous edge. Thus, there could be at most $n$ gaps, and $k$ edges from each vertex. Then, to determine if such a transformation produces an ordering not already in the ordering list, we need to search the latter. There could be up to $n^3$ orderings stored, in an index set, so using a binary search over the set is an $O(\log(n^3)) = O(\log(n))$ process, and comparing two orderings is a $O(n)$ process. There are, at most, $n^3$ such iterations in stage 3. Theoretically, stages 1--3 could be repeated $n$ times as the requirement to return from stage 3.3 to stage 1 is only that the number of gaps (at most $n$) has decreased by at least 1. Therefore the worst-case complexity of the algorithm detailed above is $O(n \times n^3 \times n \times \log(n) \times nk^4 \times n) = O(n^7 \log(n) k^4)$. For a sparse graph, where $k$ is $O(1)$, this reduces to $O(n^7 \log (n))$. Other operations (such as adding the newly constructed ordering to the ordering list) are dominated by the time taken to search for the transformations in stage 3.

Our experience indicates, however, that this worst-case complexity is very unlikely to be encountered in practice. Experimentally, we have seen that when SLH reaches stage 3 the number of gaps is very close to the minimal possible number, so stages 1--3 only need to be repeated a handful of times. Also, the likelihood of needing to open many gaps before one can be closed in stage 3 is extremely small. So, even in cases where a non-Hamiltonian graph is submitted to SLH the performance is likely to be closer to $O(n^5 \log(n) k^4)$, or $O(n^5 \log(n))$ for a sparse non-Hamiltonian graph. Furthermore, in our experiments we have seen that almost all Hamiltonian graphs are solved without ever needing to reach stage 2. In such cases the heuristic has complexity $O(n^3 + n^2 k^4)$, or $O(n^3)$ for sparse Hamiltonian graphs. This bound arises because there are at most $n^2$ iterations in stage 1, at each iteration we search from a prescribed gap so there are $O(k^4)$ choices of floating transformations, and after each iteration we need to record the new ordering which is an $O(n)$ process.

\section{Motivation}\label{sec-motivation}
It is worth mentioning that the performance of SLH has exceeded
our expectations. By that we mean that SLH correctly solved all
graph instances that we tried, which contained up to 5000 vertices. These instances included several
particularly hard graphs such as generalized Petersen (GP) cubic graphs
that contain only 3 Hamiltonian cycles, clique graphs, leap graphs, TSPLIB
graphs and many others, in an entirely deterministic
way. That is, we have not yet encountered a Hamiltonian graph where SLH fails
to find a Hamiltonian cycle before the invocation of the stopping rules.
Moreover, SLH works with the given arrangement of the initial ordering without benefiting
from randomisation techniques and preprocessing.

The effectiveness of SLH can arguably be justified by our new perspective on what it means to ``improve"
after each iteration. Let us suppose that an undirected graph $G$ is represented
by the current ordering $\mathcal{C}$ that contains $g(\mathcal{C})\ge 0$ gaps, and that a Hamiltonian
ordering $\mathcal{C}^H$ corresponding to a Hamiltonian cycle $H$ is known. Also, assume that there are $k$ snakes on $\mathcal{C}$ that
 are edges of the Hamiltonian cycle $H$ (or, equivalently, snakes in $\mathcal{C}^H$). We define the {\em distance} between $\mathcal{C}$ and $\mathcal{C}^H$ as
$$ \operatorname{dist}(\mathcal{C},\mathcal{C}^H)= n-k.$$ The {\em difference} between $\mathcal{C}$
and $\mathcal{C}^H$ is then defined as
$$ \Delta(\mathcal{C},\mathcal{C}^H)= \frac{g(\mathcal{C})}{3} +  \operatorname{dist}(\mathcal{C},\mathcal{C}^H).$$
In our perspective, a transformation that maps ordering $\mathcal{C}$ to ordering $\mathcal{C}'$ improves the tour
if $\Delta (\mathcal{C}',\mathcal{C}^H) < \Delta(\mathcal{C},\mathcal{C}^H)$ and hence, constitutes a desirable transformation.
Specifically, by applying such transformations iteratively we ensure $\mathcal{C}^H$ will eventually be obtained. We now show that from any ordering $\mathcal{C} \neq \mathcal{C}^H$
there exists a $\gamma$ or $\varkappa$ transformation, or a composition $\gamma \circ \varkappa$ that,
in our sense, improves the tour. This approach is in contrast with $k-$opt algorithms
(adapted for HCP), where the focus is solely on reducing the number of gaps. We argue that reducing the number of
gaps is not a sufficient measure of improvement, because closing gaps alone
might not bring us any closer to $\mathcal{C}^H$; indeed, the $\operatorname{dist}(\mathcal{C}',\mathcal{C}^H)$ might increase
compared to $\operatorname{dist}(\mathcal{C},\mathcal{C}^H)$. Thus, in situations where the distance is large, $k-$opt algorithms might struggle with finding a Hamiltonian cycle, even if the number of gaps is small.
In Section \ref{sec-performance} we present examples when $k-$opt approach
quickly reduces the number of gaps to 1 but fails to find any Hamiltonian cycle $H$,
arguably because the distances of any Hamiltonian ordering to the current ordering with 1 gap is quite large. By contrast, reducing the distance, without controlling the number of gaps,
 is enough to converge to a Hamiltonian ordering. However,
since we do not know the Hamiltonian cycle $H$, we can not measure the distance. Hence, reducing the number of
gaps in an ordering is merely a pragmatic, surrogate, objective. Nonetheless, SLH is willing to sacrifice the latter to move,
\lq\lq laterally", to a possibly better ordering.

We now proceed to the argument that $\gamma$ and $\varkappa$ transformations (or a composition of the two)
suffice to transform an ordering $\mathcal{C} \neq \mathcal{C}^H$ to another ordering $\mathcal{C}'$ such that $\Delta (\mathcal{C}',\mathcal{C}^H) < \Delta(\mathcal{C},\mathcal{C}^H)$.

More precisely, we show that the $\operatorname{dist}(\mathcal{C}',\mathcal{C}^H) \leq \operatorname{dist}(\mathcal{C},\mathcal{C}^H) - 1$ and
that $g(\mathcal{C}') \leq g(\mathcal{C}) + 2$.

\begin{figure}[h!]\begin{center}\includegraphics[scale=0.31]{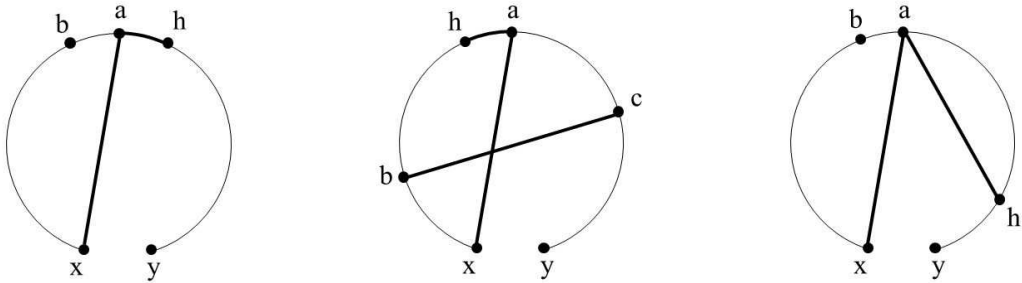}\\
\caption{Options for a Hamiltonian cycle containing edge $(x,a)$ constituting a ladder in a given ordering.}\label{fig: SLH_clockwise}\end{center}\end{figure}

Suppose that $g(\mathcal{C}) > 0$ and that $\mathcal{C}$ contains a
gap $(x,y)$. Then at least one edge in $H$ must be a ladder in $\mathcal{C}$, say $(x,a)$, emanating from $x$.
There are three options for the other edge $(a,h)$ with endpoint $a$ on the Hamiltonian cycle $H$:
\begin{itemize}
\item[(i)] $(a,h)$ is a snake in $\mathcal{C}$, and the segment $(a\smallfrown h,\ldots,y)$ does not contain $x$;

\item[(ii)] $(a,h)$ is a snake in $\mathcal{C}$, and the segment $(a\smallfrown h,\ldots,y)$ contains $x$;

\item[(iii)] $(a,h)$ is a ladder in $\mathcal{C}$.
\end{itemize}
In situations (i) and (iii) if we choose
$\gamma(y,x,a)$ it results in the ordering $\mathcal{C}'$, such that $\Delta(\mathcal{C}',\mathcal{C}^H) \le \Delta(\mathcal{C},\mathcal{C}^H)-1,$ because the distance $\operatorname{dist}(\mathcal{C}',\mathcal{C}^H) \le
\operatorname{dist}(\mathcal{C},\mathcal{C}^H) -1$. This is because the ladder $(x,a)$, which is in $H$, becomes a snake, whereas if snake $(a,b)$ exists, it becomes a ladder but is not in $H$. Also, the number of
gaps in $\mathcal{C}'$ does not increase compared to $\mathcal{C}$, though it could decrease if $(y,b)$ was a ladder, or $(a,b)$ was a gap.

In situation (ii) $H$ continues from vertex $h$ to some vertex $b$ on the segment $(a\smallfrown h,\ldots,x)$ and then moves to a vertex $c$ on the segment $(h\smallfrown a,\ldots,y)$.
Without loss of generality, let $\mathcal{C}$ be $(x,\ldots,d,b,r\ldots,h\smallfrown a,e,\ldots,f,c,\ldots,y|x)$. There are now three possibilities

\begin{figure}[h!]\begin{center}\includegraphics[scale=0.31]{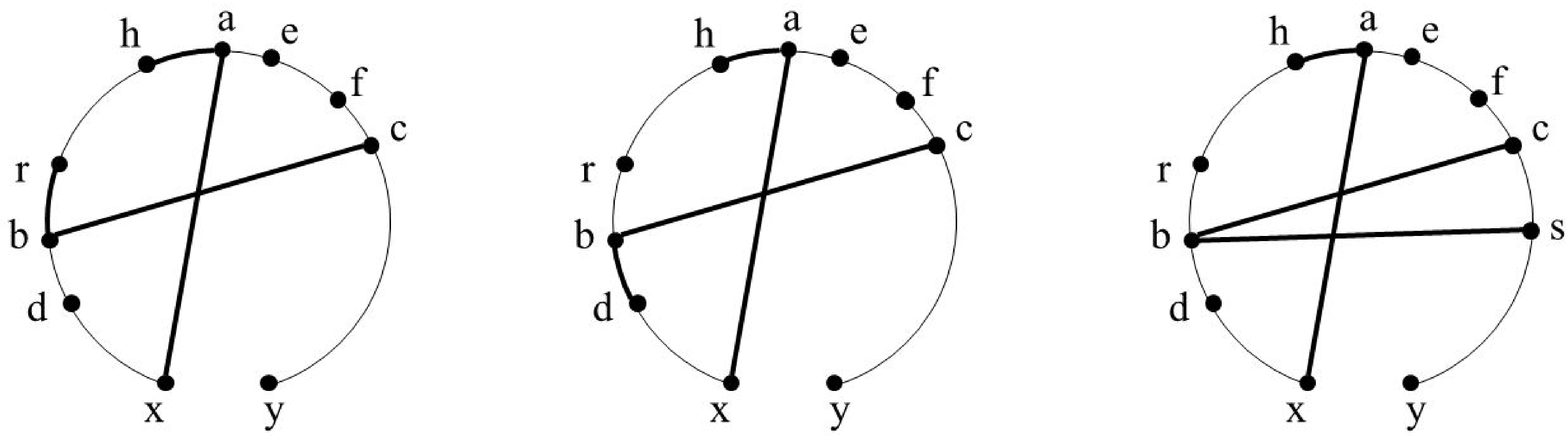}\\
\caption{From left to right, Hamiltonian cycle containing edge $(b,r)$, $(b,d)$ and $(b,s)$ respectively.}\label{fig: anticlock3}\end{center}\end{figure}

\begin{enumerate}
\item Edge $(b,r) \in H.$ In this case we apply $\varkappa(x,a,b,c)$ transform that maps $\mathcal{C}$ into some $\mathcal{C}'$. The number of gaps will not increase by more than $1$.
However, $\mathcal{C}'$ contains new snakes $(x,a)$ and $(b,c)$ which are in $H$,
but does not contain $(c,f)$ that is a snake in $\mathcal{C}$, and may be an edge in $H$.
Snake $(b,d)$ that becomes a ladder on $\mathcal{C}'$ does not belong to $H$
and therefore its absence on the circle does not increase $\Delta.$ So, overall, the number of snakes in $\mathcal{C}'$ that are in $H$ must grow by at least 1 compared to $\mathcal{C}$.

\item Edge $(b,d) \in H.$ In this case we apply $\gamma(e,d,b)\circ \varkappa(x,a,b,c)$ transform that maps $\mathcal{C}$ into  $\mathcal{C}'$. Again, the number of gaps can not increase by more than $1$.
On the other hand, $\mathcal{C}'$ contains former ladders $(x,a)$ and $(b,c)$ (which are edges in $H$) as snakes,
but former snake $(c,f)$ is a ladder in $\mathcal{C}'$, and may be an edge in $H$. Former snake $(b,r),$ that becomes
a ladder on $\mathcal{C}'$, does not belong to $H$ and therefore does not contribute to $\Delta.$ So, overall, the number of
snakes in $\mathcal{C}'$ that are in $H$ must grow by at least 1 compared to $\mathcal{C}$.

\item Some edge $(b,s) \in H$ despite being a ladder in $\mathcal{C}$. In this case $\varkappa(x,a,b,c)$ turns ladders $(x,a)$ and $(b,c)$, both of which are in $H$, into snakes in $\mathcal{C}'$. Former snake $(c,f)$ becomes a ladder
in $\mathcal{C}'$ and may be in $H$. All other former snakes that become ladders in $\mathcal{C}'$ are not in $H$. So, overall, the number of snakes in $\mathcal{C}'$ that are in $H$ must grow by at least 1 compared to $\mathcal{C}$.

\end{enumerate}

In all three cases we see that $\operatorname{dist}(\mathcal{C}',\mathcal{C}^H) \le \operatorname{dist}(\mathcal{C},\mathcal{C}^H) -1$, and $g(\mathcal{C}') \leq g(\mathcal{C}) + 2$. Therefore,
$\Delta(H,\mathcal{C}) -\Delta(H,\mathcal{C}') \ge 2/3,$ even though the transformation that maps $\mathcal{C}$ to $\mathcal{C}'$
might increase the number of gaps.

Suppose now that $g(\mathcal{C})=0$, that is $\mathcal{C}$ corresponds to a Hamiltonian cycle, different from $H.$
Then $H$ must contain an edge, say $(x,a)$, that is a ladder in $\mathcal{C}$. Then, there is vertex $y$, adjacent to $x$ in $\mathcal{C}$, such that the snake $(x,y)$ is not in $H.$ If we treat $(x,y)$ as if it were a gap (since we are not concerned if it is transformed into a ladder), and  apply exactly those transformations as above in the corresponding situation, we improve the distance by at least 1 in every case and may create at most 2 gaps if $\varkappa$ or a composition of $\varkappa$ and $\gamma$ is applied. Therefore, an overall improvement for the $\Delta$ function will be at least $1/3,$ compared to $2/3$
when $(x,y)$ is a genuine gap.

The above argument implies that, as mentioned in the introduction, for every ordering $\mathcal{C}$ and any Hamiltonian ordering $\mathcal{C}^H \neq \mathcal{C}$,
there exists a transformation $T$ mapping $\mathcal{C}$ to $\mathcal{C}^H$, such that
\begin{equation*}\label{comp_gk}
T=\prod_{i=1}^\ell\gamma^{\varepsilon_i} \circ \varkappa^{\delta_i},
\end{equation*}
where $\varepsilon_i$ and $\delta_i$ can take values $0$ or $1$, and $\ell$ bounds the number of transformations. At any stage of the process, there exists at least one transformation that can decrease $\Delta$ by at least $1/3$. Performing the correct sequence of these transformations in the process of constructing $T$ ensures convergence to $\mathcal{C}^H$ in no more than $3\Delta$ transformations.

Of course, in our implementation of the algorithm, we do not know a Hamiltonian cycle in advance, so we use the transformations outlined in
Section \ref{sec-motivation}. Each of those transformations can be represented as $\gamma^{e_i} \circ \varkappa^{d_i}$ where $e_i, d_i \in \{0, 1, 2\}$,
and we stop after $n^3$ such transformations. So, given a graph and an initial ordering, our implementation of SLH produces, in a deterministic fashion, a more specialised mapping $T$
\begin{equation*}\label{algorithm_mapping}
T=\prod_{i=1}^L\gamma^{e_i} \circ \varkappa^{d_i},
\end{equation*}
where $e_i, d_i \in \{0, 1, 2\}$ depend upon which transformation is performed in step $i$, and $L \leq n^4$. Then, either $T$ transforms the initial ordering to a Hamiltonian ordering,
or we report that the graph is \lq\lq likely non-Hamiltonian".

\begin{figure}[h]\begin{center}\includegraphics[scale=0.24]{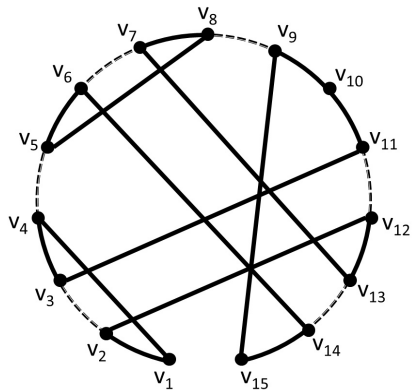}\\
\caption{A 15 vertex graph with only one Hamiltonian cycle. The dashed lines are the edges that are not a part of the Hamiltonian cycle.}\label{fig: 15node}\end{center}\end{figure}

Next, we use the 15-vertex example given in Figure \ref{fig: 15node} to demonstrate the contrast between SLH and the most famous
$k$-opt heuristic; the classical Lin-Kernighan algorithm \cite{LK}. The main tool for finding the optimal
 tour (Hamiltonian cycle in this case) in the Lin-Kernighan algorithm is using sequential exchanges where an improvement is
 possible. In other words, consider a set of snakes and gaps $X = \{x_1, x_2, ...,x_k\}$ and a set of ladders
 $Y = \{ y_1, y_2, ...,y_k  \}$ where the exchange of every $x_i$ with $y_i$ constitutes an improved tour (i.e. fewer gaps).
  A sequential exchange is when $x_i$ and $y_i$ share an endpoint, and so do $y_i$ and $x_{i+1}$, for all $i$. For the graph
   of Figure \ref{fig: 15node} with only one Hamiltonian cycle, an improvement is possible only if $X$ consists of the gap
    and all snakes that are not a part of the Hamiltonian cycle, while $Y$ consists of the $7$ ladders that are on the
    Hamiltonian cycle. This can be seen because there is only one gap, and therefore improvement implies that the number
    of gaps must be zero after a $k-$opt transformation. Therefore, a unique $7$--opt transformation exists that can improve
     the tour by exchanging $X$ with $Y$ on the ordering. Although it is possible to construct this $7$--opt transformation
     through a sequence of sequential edge-exchange moves (which individually do not improve the tour), the classical Lin-Kernighan
     algorithm would be unable to do so even if unlimited backtracking was permitted\footnote{It should be noted that
     modern implementations of Lin-Kernighan view sequential exchanges in a more refined sense, which permits them to construct the required $7$--opt transformation.}. In the words of the authors of \cite{LK}, ``the procedure is
      augmented by a limited defense in [such] situation[s]". So Lin--Kernighan's main technique to deal with situations
       that require complicated non-sequential exchanges is randomisation techniques and hoping that an easier situation arises. However,
        if the graph is large and the number of Hamiltonian cycles is small, randomisation might be of little help as demonstrated
         in the next section with generalised Petersen graphs. On the other hand, one can check that the following SLH
          transformation $T$ can unfold the graph and find the Hamiltonian cycle in the graph of Figure \ref{fig: 15node}.
\begin{equation*}
T = \gamma(v_{13},v_{9},v_{15}) \circ \gamma(v_{9},v_{6},v_{14}) \circ \varkappa(v_{5},v_{8},v_{7},v_{13}) \circ \gamma(v_{15},v_{10},v_{11}) \circ \gamma(v_{5},v_{2},v_{12}) \circ \varkappa(v_{1},v_{4},v_{3},v_{11}).
\end{equation*}
The current implementation of SLH solves this graph in stage 1. However, more difficult instances are sometimes solved
in stages 2 and 3, where opening transformations are featured. In those stages, the differences between SLH and the
HCP adapted L--K algorithm are best revealed.

\section{Experimental results}\label{sec-performance}

The times reported in this section were obtained by running SLH on a Dell R210 PC, (Intel(R) Xeon(R) E3-1270, 3.4GHz, 16GB RAM) running Linux CentOS 6.

We first tested SLH on several million subcubic
(degree three or less) graphs containing up to 50 vertices.
 These tests helped us design effective stages and transformations
  for the algorithm. With the stages in place and the transformations selected as described in Section \ref{sec-algorithm},
  the performance of SLH proved to be competitive with benchmark solvers, based upon experiments conducted on a variety of graphs. Specifically, we compared
 the performance of SLH to that of {\em Concorde} TSP Solver \cite{Concorde},
 Helsgaun's Lin-Kernighan implementation, {\em LKH} HCP solver \cite{Helsgaun},
  and Eppstein's HCP deterministic solver \cite{eppstein} for cubic graphs. Though Concorde is
  a TSP solver, we represent instances of HCP as TSP by assigning weight $0$ to existing edges
  (edges of the given graph) and weight $1$ to non-existing edges. It should be noted that LKH has a built-in pre-solve phase that performs subgradient optimisation. In order to obtain a fair comparison with SLH, we have chosen to evaluate LKH both with and without this pre-solve stage. As is demonstrated below, the pre-solve is very effective on certain classes of graphs, which indicates that SLH would benefit from a similar pre-solve phase. Some experimentation with alternative parameters was performed, however with only one exception, we did not discover any instances where altering the parameters significantly improved the solving time or reliability in any of the comparison solvers. The sole exception was in switching off the restricted search option in LKH; in this case, LKH was more reliable, but took longer to solve. For the remainder of this manuscript, default parameters are assumed. We ran the Windows/Cygwin implementation of Concorde, and Eppstein's solver, on a Dell Optiplex 780 (Intel(R) Core(TM)2 Quad CPU Q9650, 3.0Ghz, 4GB RAM) running Windows 7. We ran LKH v2.0.7 on the same Dell R210 PC used to run SLH.

For the test graphs we selected all TSPLIB HCP data \cite{TSPLIB}; Flower snarks \cite{flower} ${\bf J}_k$ for $k = 5, 15, 25, 35$; Sheehan graphs \cite{sheehan} for $n = 50, 60, 70, 80$; and also generalised Petersen cubic graphs \cite{GP} $GP(n,2)$ for $n=39,45,51,63,123,243.$ The TSPLIB graphs are benchmark instances notable for their large size (up to 5000 vertices). They are all sparse, irregular graphs, and are all Hamiltonian. The Flower snarks are a family of non-Hamiltonian cubic graphs of size $n = 4k$, that all contain Hamiltonian paths; so the optimal solution, from a TSP viewpoint, is a path of length 1. These were tested to give an indication of how quickly SLH can detect non-Hamiltonicity compared to the other solvers tested, which are not designed to deal with such instances. This is an important test, as the essential difficulty of HCP is that a Hamiltonian cycle may not exist for a given graph, and hence it is interesting to see how the algorithms will perform in these cases. Sheehan graphs are maximally dense uniquely Hamiltonian graphs; that is, Hamiltonian graphs with only a single Hamiltonian cycle, that contain as large a ratio of edges to vertices as possible. Finally, the generalised Petersen graphs, for the values of $n$ we chose (that is, $n \equiv 3 \mod 6$), are cubic graphs that contain only three Hamiltonian cycles, which is the minimum number a Hamiltonian cubic graph may contain. For the Sheehan graphs and the generalised Petersen graphs we chose a random permutation of the vertices to help disguise the inherent structure.

In our tests, the success rate of LKH algorithm -- which invokes randomisations -- is determined after many trials (50 trials for Sheehan graphs, 1000 trials for other graphs) and execution time is taken as average over the trials. Eppstein's algorithm is only designed for cubic graphs and so was tested only on cubic graphs, that is, the Flower snarks and the generalised Petersen graphs. Eppstein's algorithm finds all Hamiltonian cycles, so in the case of Hamiltonian graphs, the total time was divided by the number of Hamiltonian cycles. Also note that Eppstein's algorithm, like SLH (and unlike Concorde and LKH), is a deterministic algorithm.

The results of the above comparisons are reported in Table \ref{tab: performance}. Assessment {\em Fail} means that the solver terminated without managing to find a Hamiltonian cycle despite the graph being Hamiltonian, while {\em Time Fail} means that the solver had not concluded within 48 hours. Note that, for the non-Hamiltonian instances (the Flower snarks), assessment {\em Fail} does not apply as there is no Hamiltonian cycle to find.

\begin{table}[h]
\hspace*{-0.6cm}
\begin{tabular}{|l||c|cc|cc|cc|c|}
\hline
\multirow{2}{*}{{\bf Graphs}} & {\bf Concorde} & \multicolumn{2}{c|}{{\bf LKH}} & \multicolumn{2}{c|}{{\bf LKH without presolve}} & \multicolumn{2}{c|}{{\bf SLH}} & {\bf Eppstein} \\
& Time(sec) & Time(sec) & Success (\%) & Time(sec) & Success (\%) & Time(sec) & Stage & Time(sec)\\  \hline
ALB 1000  & 4.95       & 0.0  & 100 & 0.0  & 100 & 0.16   & 1 & - \\
ALB 2000  & 7.30       & 0.0  & 100 & 0.0  & 100 & 0.79   & 1 & - \\
ALB 3000a & 9.56       & 0.1  & 100 & 0.1  & 100 & 3.19   & 1 & - \\
ALB 3000b & 9.94       & 0.1  & 100 & 0.1  & 100 & 3.31   & 1 & - \\
ALB 3000c & 9.95       & 0.1  & 100 & 0.1  & 100 & 2.91   & 1 & - \\
ALB 3000d & 10.14      & 0.1  & 100 & 0.1  & 100 & 3.26   & 1 & - \\
ALB 3000e & 10.44      & 0.1  & 100 & 0.1  & 100 & 2.77   & 1 & - \\
ALB 4000  & 13.45      & 0.1  & 100 & 0.1  & 100 & 5.75   & 1 & - \\
ALB 5000  & 17.24      & 0.1  & 100 & 0.1  & 100 & 12.48  & 1 & - \\ \hline
Flower5   & 0.15       & 0.0  & 100 & 0.0  & 100 & 0.04   & 3 & 0.57 \\
Flower15  & 2014.01    & 0.1  & 100 & 0.1  & 100 & 1.80   & 3 & 14.31 \\
Flower25  & Time Fail  & 0.2  & 100 & 0.2  & 100 & 12.26  & 3 & 11675.73 \\
Flower35  & Time Fail  & 0.3  & 100 & 0.3  & 100 & 41.91  & 3 & Time Fail \\ \hline
Sheehan50 & 0.23       & 0.0  & 100 & 2.8 & 100 & 0.02   & 1 & - \\
Sheehan60 & 0.33       & 0.0  & 100 & 16.5 & 100  & 0.27   & 1 & - \\
Sheehan70 & 0.41       & 0.0  & 100 & 69.2 & 96 & 0.29   & 1 & - \\
Sheehan80 & 0.59       & 0.0  & 100 & & 24  & 1.79   & 1 & - \\ \hline
GP(39,2)  & 34.38      & 0.0  & 100  & 0.0  & 100 & 0.01  & 1 & 0.6 \\
GP(45,2)  & 50.91      & 0.0  & 100  & 0.0  & 100 & 0.06  & 1 & 3.4 \\
GP(51,2)  & 50.09      & 0.0  & 99.8 & 0.0  & 99.8 & 0.07  & 1 & 9.5  \\
GP(63,2)  & 737.88     & 0.0  & 84.9 & 0.0  & 83.8 & 0.15  & 1 & 129.1  \\
GP(123,2) & Fail       & 0.3 & 00.1 & 0.3 & 00.1 & 6.27  & 2 & Time Fail \\
GP(243,2) & Fail       & 0.9 & 00.0 & 0.9 & 00.0 & 620.48  & 3 & Time Fail \\ \hline
\end{tabular}
\caption{Comparative performance of SLH for many large or difficult instances of HCP.}
\label{tab: performance}
\end{table}

Table \ref{tab: performance} reveals that the performance of SLH compares increasingly favourably as the difficulty of the graph grows. For general instances, LKH dominates both Concorde and SLH, but struggles to solve reliably and efficiently when there is only a small number of orderings with the minimal number of gaps. Although LKH is able to solve the Sheehan graphs very quickly, it does so in its pre-solve phase. As Table \ref{tab: performance} demonstrates, its performance on these graphs deteriorates when this phase is disabled. Of course, since SLH does not currently have any pre-solve phase, arguably, it is more meaningful to compare it to LKH without pre-solve. The latter, like SLH, is essentially a $k$-exchange type algorithm. For the Flower snarks, which are hypohamiltonian, Concorde quickly finds a Hamiltonian path (the optimal solution) but takes much longer to determine that no Hamiltonian cycle exists.

For the two largest tested instances of the generalised Petersen graphs, SLH had to proceed to later stages, where opening transformations were performed, in order to find a Hamiltonian cycle. For these two instances, Eppstein's algorithm was unable to provide a solution within 48 hours, and in both cases Concorde encountered a numerical error that prevented it from returning any result at all. LKH was able to find a solution to GP(123,2), but only in nine out of 1000 attempts (or eleven attempts without the pre-solve phase), and failed to produce any solutions to GP(243,2) in 1000 attempts, with or without the pre-solve phase.

For the cubic graphs, we ran Eppstein's algorithm, which was dominated by the performance of the other three algorithms; however, it should be noted that Eppstein's algorithm enumerates {\em all} Hamiltonian cycles, rather than terminating once one is found. The adaptation of SLH to find additional Hamiltonian cycles, after a first has been identified, is a topic for future research.

For all four solvers tested, the hardest graphs in our experiments to solve were the {\em generalised Petersen} cubic graphs $GP(n,2)$ for $n \equiv 3 \mod 6$. It is interesting
       to note that, although LKH failed to find a Hamiltonian cycle in $GP(243,2)$,
       it very quickly descended to an ordering $\mathcal{C}$ that contained exactly
        one gap. At the same time, the value of $\Delta(\mathcal{C},\mathcal{C}^H)$ remained large for all
        three possible choices of $\mathcal{C}^H$ due to the small number of common snakes in $\mathcal{C}$ and any $\mathcal{C}^H$. In our experiments,
        Helsgaun's LKH generated a final ordering $\mathcal{C}$ for $GP(243,2)$, for which $\Delta(\mathcal{C},\mathcal{C}^H)$
        was $150\frac{1}{3}$, $154\frac{1}{3}$ and $190\frac{1}{3}$ for the three Hamiltonian cycles.

Table \ref{tab: success} reports the success rate of SLH over large test sets. SLH was first tested on the Foster census (see Bouwer et al. \cite{bouwer} and Royle \cite{Royle}). SLH was next tested on 1,000,000 cubic graphs of size 100 (generated uniformly, at random, see \cite{wormald}), separated into four test sets. Finally, SLH was tested on 10,000 cubic graphs of size 1,000 (generated uniformly by random). The success rate of SLH was checked by confirming that the reported Hamiltonian cycles exist in the graph\footnote{Of course, this is guaranteed by the algorithm, but was checked nonetheless to confirm our implementation was accurate.}, and independently confirming any result of non-Hamiltonicity. There were 95 non-Hamiltonian graphs in the 1,000,000 randomly generated cubic graphs of size 100. All of the randomly generated cubic graphs of size 1,000 were Hamiltonian. The reported times are the total times to solve all instances in each set of graphs, with the input/output times removed.

\begin{table}[h]
\centering
\begin{tabular}{|l||c|c|c|c|} \hline
{\bf Graphs} & {\bf Number} & {\bf Size range} & {\bf Total time (s)} & {\bf Success rate}\\ \hline
Foster Census & 332 & 4 -- 998 & 17.74 & 100\\ \hline
Cubic100 \#1 & 250,000 & 100 & 448.39 & 100\\
Cubic100 \#2 & 250,000 & 100 & 439.78 & 100\\
Cubic100 \#3 & 250,000 & 100 & 476.10 & 100\\
Cubic100 \#4 & 250,000 & 100 & 440.40 & 100\\\hline
Cubic10000   & 10,000  & 1000 & 31.72 & 100\\\hline
\end{tabular}
\caption{Performance and success rate of SLH for large sets of graphs.}
\label{tab: success}
\end{table}

As can be seen in Table \ref{tab: success}, SLH succeeded in finding a Hamiltonian cycle in all Hamiltonian graphs tested. To date the authors have not found an instance of a Hamiltonian graph for which SLH has terminated without finding a Hamiltonian cycle, having tested graphs of sizes up to 5000.

\section{Conclusions} \label{sec-conclusion}

In this paper we presented SLH; a new heuristic with polynomial complexity for solving the Hamiltonian Cycle Problem. Table \ref{tab: performance} shows that SLH, while still in its infancy, is already a competitive HCP solver, while Table \ref{tab: success} indicates that SLH is extremely reliable, succeeding in finding a Hamiltonian cycle in every single tested Hamiltonian graph, notably including graphs possessing very few Hamiltonian cycles.  This includes graphs where other, benchmark, solvers failed to find any such cycles. Furthermore, these results were achieved without the need for randomisation techniques that are a common feature in many contemporary HCP heuristics. We believe that balancing  the use of our $\gamma$ and $\varkappa$ generator transformations overcomes some of the difficulties that other heuristics try to overcome by randomisation.

In addition, we expect that in future versions, the speed of SLH will be significantly improved by optimising the implementation in a variety of ways.  In particular, replacing SLH-closing transformations with
Lin-Kernighan type sequential transformations, such as those implemented in LKH \cite{Helsgaun}, seems promising and deserves further exploration. Similarly, investigations of both alternative initial orderings and alternative stopping rules could lead to improvements in performance. Taking advantage of a pre-solver, such as the subgradient optimisation routine used in LKH, is also likely to improve the solving time of the algorithm in general. A demonstration version of SLH is available on our website \cite{SLH}, where users are invited to submit and solve their HCP problems through our online interface.

\section*{Acknowledgments}

The authors gratefully acknowledge useful comments from the anonymous referees which improved the exposition, and useful discussions with Brendan McKay and Gordon Royle that helped us to find suitable test instances. The editor, William Cook, also contributed significantly by suggesting further testing and changes of inaccurate statements. The research presented in this manuscript was supported by the ARC Discovery Grant DP120100532.

\bibliographystyle{plain}   

\begin{thebibliography}{99}
\bibitem{Concorde} Applegate, D.L., Bixby, R.B., Chav\'{a}tal, V., and Cook, W.J. Concorde TSP Solver: http://www.tsp.gatech.edu/concorde/index.html.
\bibitem{ConcordeBook} Applegate, D.L., Bixby, R.B., Chav\'{a}tal, V., and Cook, W.J. {\em The Traveling Salesman Problem: A Computational Study}. Princeton University Press (2006).
\bibitem{SLH} Baniasadi, P., Clancy, K., Ejov, V., Filar, J.A., Haythorpe, M., and Rossomakhine, S. Snakes and Ladders Heuristic -- Web Interface: http://fhcp.edu.au/slhweb/ (2012).
\bibitem{bouwer} Bouwer, I.Z., Chernoff, W.W., Monson, B., Star, Z. {\em The Foster Census}. Charles Babbage Research Center, Winnipeg (1988).
\bibitem{eppstein} Eppstein, D. The Traveling Salesman Problem for Cubic Graphs. In Frank Dehne, J\"{o}rg-R\"{u}diger Sack, and Michiel Smid, editors, {\em Algorithms and Data Struct.}, volume 2748 of {\em Lecture Notes in Computer Science}, pages 307--318. Springer Berlin (2003).
\bibitem{Flood} Flood, M.M. The Traveling Salesman Problem. {\em Oper. Res.} 4, 61--75 (1956).
\bibitem{Gutin} Gutin, G. and Punnen, A.P. {\em Traveling Salesman Problem and Its Variations}. Kluwer Academic Publishers (2002).
\bibitem{Helsgaun} Helsgaun, K. An Effective Implementation of Lin-Kernighan Traveling Salesman Heuristic. {\em Eur. J. Oper. Res.} 126, 106--130 (2000).
\bibitem{flower} Isaacs, R. Infinite Families of Nontrivial Trivalent Graphs Which Are Not Tait Colorable. {\em Amer. Math. Monthly} 82:221--239 (1975).
\bibitem{Lawl85Traveling} Lawler, E.L., Lenstra, J.K., Rinooy Kan A.H.G., and Shmoys, D.B. {\em The Traveling Salesman Problem: A Guided Tour of Combinatorial Optimization}. John Wiley and Sons (1985).
\bibitem{Lin} Lin, S. Computer Solutions of The Traveling Salesman Problem. {\em The Bell Systems Tech. J.} 44, 2245--2269 (1965).
\bibitem{LK} Lin, S. and Kernighan, B.W. An Effective Heuristic Algorithm for The Traveling Salesman Problem. {\em Oper. Res.} 21, 496--516 (1973).
\bibitem{Royle} Royle, G., Conder, M., McKay, B., and Dobscanyi, P. Cubic symmetric graphs (The Foster Census): http://school.maths.uwa.edu.au/$\sim$gordon/remote/foster (2001).
\bibitem{sheehan} Sheehan, J. Graphs with exactly one hamiltonian circuit. {\em J. Graph Th.} 1:37--43 (1977).
\bibitem{TSPLIB} TSPLIB. Hamiltonian cycle problem (HCP): http://comopt.ifi.uni-heidelberg.de/software/TSPLIB95 (2008).
\bibitem{GP} Weisstein, E.W. Generalized Petersen Graph (From MathWorld -- A Wolfram Web Resource): http://mathworld.wolfram.com/generalizedpetersengraph.html.
\bibitem{wormald} Wormald, N. Models of Random Regular Graphs, in {\em Surveys in Combinatorics}, Cambridge University press, pages 239--298 (1999).

\end{thebibliography}

\end{document}